\magnification\magstep1
\baselineskip = 18pt
\overfullrule = 0pt
\def\n{\noindent}
\def\qed{{\hfill{\vrule height7pt width7pt
depth0pt}\par\bigskip}}

\def\n{\noindent}
\def\ms{\medskip}

\def\RR{{\mathop{{\rm I}\kern-.2em{\rm R}}\nolimits}}
\def\TT{{\mathop{{\rm I}\kern-.4em{\rm T}}\nolimits}}
\def\FF{{\mathop{{\rm I}\kern-.2em{\rm F}}\nolimits}}
\def\NN{{\mathop{{\rm I}\kern-.2em{\rm N}}\nolimits}}
\def\EE{{\mathop{{\rm I}\kern-.2em{\rm E}}\nolimits}}
\def\CC{{\rm C\kern-.18cm\vrule width.6pt height 6pt depth-.2pt
\kern.18cm}}
\def\pf{ \medskip \n {\smcaps Proof.~~}}

\font\smcaps=cmcsc10

\centerline{\bf ON A QUESTION OF NIELS GR\O NBAEK} \bigskip

\centerline {By {\smcaps Gilles Pisier}\footnote*{Partially supported by
the NSF}}
\centerline{ Texas A\&M University, USA}
\centerline{and}
\centerline{Universit\'e Paris VI, Paris}

\ms \ms \ms

\parindent=0pt 

\centerline{\smcaps Abstract}

  Let $F(X)$ denote the norm closure of the space of all finite rank
operators on a Banach space $X$. 
We show that there are Banach spaces $X$
for which the product map $a\otimes b\to ab$ does not define a surjective map from the
projective tensor product $F(X) \widehat\otimes F(X)$ onto $F(X)$.

  \vfill\eject

\baselineskip = 18pt

Let $X$ be a Banach space and let $B(X)$ be the space
of all bounded operators on $X$
equipped with the usual norm. Let $F(X)$ denote the closure in $B(X)$ of
the set of  all finite rank maps on $X$.

Niels Gr\o nbaek asked me the following question:\ Is it true that every Banach 
space has the following property
{\parindent=20pt
\item{{\bf (P)}} The product map $F(X) \widehat\otimes F(X)\to F(X)$ is onto, or 
equivalently there exists a constant $C$ such that for any finite rank map 
$u\colon \ X\to X$, there are finite rank maps $v_n\colon \ X\to X$ and 
$w_n\colon 
\ X\to X$ such that $u = \sum^\infty_1 v_nw_n$ and $\sum^\infty_1 \|v_n\|\, 
\|w_n\| \le C\|u\|$.\smallskip}

As observed by Gr\o nbaek, any space with the bounded approximation property
has this property, while, at the other end of the spectrum,
any Banach space $X$ such that $X\widehat\otimes X=X\check\otimes X$ (as constructed in [4])
also does satisfy (P) but fails the approximation property.
However,  as shown by the next statement, 
which is the main result of this note, there are intermediate cases
for which (P) fails. This  seems to be the first application
of the variant of the construction in [4] for the case of cotype $q$ with $q\not=2$.

\proclaim Theorem 1. There is a (separable) Banach space failing the above property 
(P). More precisely, for any $q>2$ any Banach space $X_q$ which is of 
cotype $q$ and such that the dual $(X_q)^*$ is a $GT$ space of cotype 2 (see below) but which
contains $\ell_q$ isomorphically, must fail (P).
\medskip

Let us briefly recall some terminology.
 Let $u\colon\ X\to Y$ be an operator between Banach spaces. Let $0<p\le q<\infty$.
Then $u$ is called  $(q,p)$-absolutely summing
(resp. $p$-absolutely summing in case $q=p$),
 if there is
a constant $C$ such that for any finite sequence $(x_i)$ in $X$, we have
$$(\sum \| u(x_i)\|^q)^{1/q} \le \sup\{ (\sum |\xi(x_i|^p)^{1/p} \mid \xi \in X^*, \|\xi\|\le
1\}.$$
We denote by $\pi_{q,p}(u)$ (resp. $\pi_{p}(u)$ in case $q=p$) the smallest constant $C$ for
which this holds. Actually, for short we will use the terms
 $(q,p)$ summing  or $p$-summing.

\n A Banach space $X$ is called of cotype $q$ ($q\ge2$) if there is a constant
 $C$ such that for any finite sequence $(x_i)$ in $X$, we have
$$\left(\sum \|   x_i \|^q\right)^{1/q} \le \left(\int \|\sum r_i(t) x_i\|^2 dt\right)^{1/2},$$
where we have denoted by $(r_i)$ the sequence of the Rademacher functions
on the Lebesgue interval.

\n As Grothendieck's fundamental theorem shows
(see [5]) there are Banach spaces $X$
(for instance $X=L_1$) such that every bounded operator  $u\colon\ X\to \ell_2$
is automatically $1$-summing. As in [4], we call these GT-spaces.
We refer the reader to [6, 5] for background on $p$-summing,
or  $(q,p)$-summing
 operators, and cotype of Banach spaces

\medskip

\n {\it Remark 2.} Let $X$ be a Banach space. 
It is proved in [4, Prop. 1.11] that $X^*$ is 
a $GT$ space of cotype 2 iff $X$ satisfies the following:\ there is a constant 
$C$ such that, if we denote by $R\subset L_1$ the closed span of the Rademacher 
functions in $L_1$ over the Lebesgue interval, then every $v\colon \ R\to X$ admits an extension
$\tilde v\colon
\  L_1\to X$ such that $\|\tilde v\|\le C\|v\|$.

\proclaim Lemma 3. Let $X$ be of cotype $q\ge 2$, and such that $X^*$ is a $GT$ 
space of cotype 2. Then there is a constant $K'$ such that any finite rank map 
$u\colon \ X\to X$ satisfies
$$\pi_{q,2}(u) \le K'\|u\|.\leqno (1)$$

\n {\smcaps Proof.}   We refer to [3] for the
definitions and the main properties of the 
$K$-convexity constant of an operator $u$, which we denote, as in [3], by 
$K(u)$ (precisely, this is defined as the best constant
in (3) below). It is proved in [3] that if $X^*$ has cotype 2 and
$X$ has cotype
$q$  then there is a constant $K''$ such that every finite rank map $u\colon \ X\to 
X$ satisfies
$$K(u) \le K''\|u\|.\leqno (2)$$
Moreover, since $X^*$ is a $GT$ space of cotype 2, by Remark~2, we have the 
following property:\ There is a constant $C'$ such that for any $ n$  and any
$x_1,\ldots, x_n$ in $X$ there is a function
 $\Phi\in L^\infty([0,1]; X)$ such that
for all $i=1,\ldots, n$ we have
$$\int  r_i(t) \Phi \ dt = x_i\quad \hbox{and}\quad \|\Phi\|_{L^\infty(X)}\le C' 
 \max\left\{\left(\sum |\xi(x_i)|^2\right)^{1/2}\ \big|\ \xi\in 
B_{X^*}\right\}.$$ 
Using the very definition of $K(u)$, this implies
$$\eqalign{(3)    \quad
\qquad\qquad\qquad\left\|\sum r_i(t)u(x_i)\right\|_{L^2(X)} 
&\le K(u) 
\|\Phi\|_{L^2(X)}\qquad\qquad\qquad\qquad\qquad\qquad\qquad\qquad\qquad\cr
&\le K(u)C' \max_{\xi\in B_{X^*}} \left(\sum |\xi(x_i)|^2\right)^{1/2}} $$
and since $X$ is of cotype $q$, we conclude that $\pi_{q,2}(u) \le C'C''K(u)$ 
for some constant $C''$ (equal to the cotype $q$ constant of $X$).
Thus, by (2) we obtain (1).\qed

\proclaim Lemma 4. Let $q>2$. Any space $X$ of cotype $q$ such that $X^*$ is a 
$GT$ space of cotype 2, but which contains 
 a subspace isomorphic to $\ell_q$    (or which 
contains $\ell_q^n$'s uniformly, in the sense of e.g. [5, p.39])   fails the property (P)
introduced above.

\pf By a result due to K\"onig-Retherford-Tomczak (see [6, \S22] or [2]),
if $p<2$, $k\ge 2$  and $1/p= 1/q_1+...+1/q_k$
 where $2\le q_i<\infty$ for each $i=1,...,k$,
there
is a constant $B$ such that for any $k$-tuple 
$u_1,...,u_k$ of operators in $B(X)$,
we have 
$$\pi_2(u_1u_2...u_k)\le B \prod_{i=1}^k \pi_{2,q_i}(u_i).$$
Now let us assume for simplicity that $2<q<4$. 
Then if we take $k=2$ and  $q_1=q_2=q$, we have
 $p<2$ so that
this implies using (1)
$$\pi_2(vw) \le BK'^2 \|v\|\, \|w\| .$$
So if the property (P) appearing in Theorem~1 held, we would have, if $u = \sum 
v_nw_n$
$$\eqalign{\pi_2\left(\sum v_nw_n\right)  &\le BK'^2 \sum \|v_n\|\, \|w_n\|\cr
&\le CBK'^2\|u\|}$$
hence  for all $u\colon \     X\to X$ with finite rank, $\pi_2(u) \le CBK'^2\|u\| $. Now applying the 
property again we obtain by a classical
 composition result of Pietsch (see [6, p. 55])
$$\eqalignno{\pi_1\left(\sum v_nw_n\right) &\le \sum \pi_2(v_n)\pi_2(w_n)\cr
&\le (CBK'{}^2)^2 \sum\|v_n\|\, \|w_n\|\cr
&\le C(CBK'^2)^2 \|u\|,}$$
  hence we find a constant $K''$  such 
that any finite rank $u\colon\ X\to X$ satisfies
$$\pi_1(u) \le K''\|u\|,$$
and a fortiori satisfies $$\pi_1(u) \le K''\pi_2(u).$$
But now this fails when $X\supset \ell_q$, $q>2$ because, by a standard
argument (for details see the equivalence of $(iv)_a$ and $(iv)_b$ in [1, Prop. 2.1]),
 this implies
that $$\Pi_2(\ell_\infty, X)= B(\ell_\infty, X)$$
  and this obviously fails 
if the space $\ell_q$, $q>2$ embeds into $X$. Indeed,
it would immediately follow
that $\Pi_2(\ell_\infty,  \ell_q)= B(\ell_\infty,  \ell_q)$,
 but it is easy to check that the
diagonal multiplication operators from 
$\ell_\infty$ into $ \ell_q$ are in general {\it not}
$2$-summing when  $q>2$.

Similarly if $q>4$ we choose an even integer 
$m$ such that ${m\over q} <1$. By property (P), 
any finite rank $u\colon \ X\to 
X$ can be written as 
$$\sum_i v_1(i)v_2(i)\ldots v_m(i)$$
with finite rank maps $v_j(i)$ such that
$$\sum_i\|v_1(i)\|\ldots \|v_m(i)\| \le C^{m-1}\|u\|.$$
Then the preceding argument leads again to $\pi_1(u) \le K''\|u\|$ for some 
constant $K''$, which is impossible. \qed\medskip

\n {\smcaps Proof of Theorem 1.} By [4, Th. 2.4], 
any space $Y$ of cotype $q$ (resp.\ which 
is separable) can be embedded in a space $X$ 
of cotype $q$ (resp.\ separable) 
such that $X^*$ is a $GT$-space of cotype 2.
 Then, by Lemma~4 it suffices to 
take $Y=\ell_q$ to conclude.\qed

\n {\it Remark.} Let $K(X)$ denote the space
of all compact operators on
$X$.  It is apparently not   known 
whether there is a Banach space $X$ such that the natural product map
 $$K(X) \widehat\otimes K(X)\to K(X)$$ is not onto.

\vskip24pt
 \centerline{\smcaps References}
\vskip24pt

\item{[1]} E. Dubinsky, A. Pe{\l}czynski and H.P. Rosenthal,
 On Banach spaces such that
$\Pi_2({\cal L}_\infty, X)= B({\cal L}_\infty, X)$,
 {\it Studia Mathematics\/} {\bf 44} (1972), 617--48.\vskip 12pt
 
 \item{[2]} H. K\"onig, {\it Eigenvalue distribution of compact operators\/},
 Birkhauser, Basel, 1989.\vskip 12pt

\item{[3]} G. Pisier, On the duality between type and cotype, {\it 
Martingale theory in harmonic analysis and Banach spaces,
Proceedings}, Springer Lecture Notes in Mathematics, 939, Cleveland, 1981,
pp.~131--44. \vskip 12pt

\item{[4]} G. Pisier,  Counterexamples  to a conjecture of
Grothendieck, {\it  Acta Mathematics} {\bf 151} (1983), 181--209.
\vskip 12pt

\item{[5]} G. Pisier, {\it Factorization of linear operators and the
geometry of Banach spaces, CBMS}, Regional
conference of the A.M.S. {\bf 60} (1986), reprinted
with corrections 1987. \vskip 12pt

\item{[6]} N. Tomczak-Jaegermann, {\it Banach-Mazur distances
and finite-dimensional operator ideals}, Pitman-Longman,
Wiley, New York, 1989.\vskip 12pt
\bye